\documentclass{amsart}

\usepackage{amssymb, ifthen, natbib}

\newtheorem{theorem}{Theorem}[section]
\newtheorem{lemma}[theorem]{Lemma}
\newtheorem{proposition}{Proposition}[section]
\newtheorem{corollary}[theorem]{Corollary}

\theoremstyle{definition}

\newtheorem{example}[theorem]{Example}

\theoremstyle{remark}

\numberwithin{equation}{section}

\newcommand{\theor}[3][]{\begin{theorem}[#1] #3 \label{thm:#2} \end{theorem}}
\def\thm#1{{t}heorem \ref{thm:#1}}

\def\lema#1#2{\begin{lemma} #2 \label{lem:#1} \end{lemma}}
\def\lem#1{{l}emma \ref{lem:#1}}

\def\corol#1#2{\begin{corollary} #2 \label{cor:#1} \end{corollary}}
\def\cor#1{{c}orollary \ref{cor:#1}}

\def\propo#1#2{\begin{proposition} #2 \label{prp:#1} \end{proposition}}
\def\prp#1{{p}roposition \ref{prp:#1}}

\def\display#1#2{\begin{equation} #2 \label{eqn:#1} \end{equation}}
\def\eqn#1{(\ref{eqn:#1})}

\def\Section#1#2{\section{#2\label{sec:#1}}}

\def\sec#1{{s}ection \ref{sec:#1}}

\newcommand{\from}{\colon}
\newcommand{\Nat}{\mathbb{N}}

\author[\relax]{Edward J. Green}
\address{Department of Economics, The Pennsylvania State University,
  University Park, PA 16802, USA}
\email{eug2@psu.edu}

\newcommand{\hcf}{the Human Capital
  Foundation, (\texttt{http://www.hcfoundation.ru/en/}) and
  particularly Andrey P.\ Vavilov, for research support through the
  Center for the Study of Auctions, Procurements, and Competition
  Policy (CAPCP, \texttt{http://capcp.psu.edu/}) at the Pennsylvania
  State University.\ }

\newcommand{\gobble}[1]{} 
\newcommand{\Junk}{}
\newcommand{\junk}[1]{\renewcommand{\Junk}{\expandafter\gobble #1}}
\newcommand{\tc}[1]{\junk{#1} \ifthenelse{\equal{\Junk}{}}{#1^+}{(#1)^+}}

\newcommand{\rst}{\restriction}
\newcommand{\transpose}[1][{12}]{t_{#1}} 
\newcommand{\tpos}[2][{12}]{t_{#1}(#2)} 
\newcommand{\Tpos}[1]{\widetilde{#1}}
\newcommand{\cmps}{\relax} 
\newcommand{\BD}{\boldsymbol{\Delta}}
\newcommand{\BO}{\boldsymbol{\Omega}}
\newcommand{\BP}{\boldsymbol{\Pi}}
\newcommand{\BS}{\boldsymbol{\Sigma}}
\newcommand{\Part}{\mathcal{P}}
\newcommand{\N}{\mathcal{N}} 
\newcommand{\id}{i} 

\begin{document}

\title[\relax]{Embedding an
  Analytic Equivalence Relation in the Transitive
  Closure of a Borel Relation} 

\thanks{The author thanks \hcf\ Also he gratefully acknowledges the
  hospitality of CEMFI (\texttt{http://www.cemfi.es/}) where the
  research was conducted.}


\date{2012.03.07}

\begin{abstract}
The transitive closure of a reflexive, symmetric, analytic relation
is an analytic equivalence relation. Does some smaller class contain
the transitive closure of every reflexive, symmetric, closed relation?
An essentially negative answer is provided here. Every analytic
equivalence relation on an arbitrary Polish space is Borel embeddable
in the transitive closure of the union of two smooth Borel equivalence
relations on that space. In the case of the Baire space, the two
smooth relations are closed and the embedding is homeomorphic.
\end{abstract}

\maketitle

\Section a{Introduction}

This note answers a question in descriptive set theory that arises in
the context of the Bayesian theory of decisions and games. It concerns
the notion of common knowledge, formalized by Robert Aumann
\citeyearpar{Aumann-1976}. For an event $A$ that is represented as a
subset of a measurable space $\Omega$, Aumann defines the event that
an agent \emph{knows} $A$ to be the event $A \setminus [\Omega
  \setminus A]_\Part$, where $\Part$ is the agent's \emph{information
  partition} of $\Omega$.\footnote{$[A]_\Part$ denotes $\bigcup \{ \pi
  \mid \pi \in \Part \text{ and } \pi \cap A \neq \emptyset \}$, the
  saturation of $A$ with respect to $\Part$. If $E$ is an equivalence
  relation, then $[A]_E$ denotes the saturation of $A$ with respect to
  the partition induced by $E$. Aumann's definition
  corresponds to the truth condition for $\square A$ in
  \citet{Kripke-1959}.}  If $\Part$ is the meet of individual agents'
information partitions (in the lattice of partitions where $\Part' \le
\Part'' \iff \Part'' \text{ refines } \Part'$), then Aumann defines 
\display o{A \setminus [\Omega \setminus A]_\Part}
to be the event that $A$ is common knowledge among the
agents.\footnote{Aumann sketches an argument---reminiscent of a general
  principle in proof theory (cf.~\citet[Lemma 6.4.8,
    p.~89]{Pohlers-2009})---that this definition is equivalent to the
  intuitive, recursive definition of common knowledge: that $A$ has
  occurred and that, for all $n \in \Nat$, both agents know\dots
  that both agents know ($n$ times) that $A$ has occurred.}

Aumann restricts attention to the case that $\Omega$ is countable (or
that the Borel $\sigma\/$-algebra on $\Omega$ is generated by the
elements of a countable partition), so that measurability issues do
not arise. But, otherwise, the passage from information partitions to a
common-knowledge partition is very badly behaved, as is the passage
from an information partition $\Part$ and an event $A$ to the related
event that $A$ is known according to $\Part$. For example, let $X$ be an
arbitrary subset of $(0,1)$, and let $\Omega = [0,2]$. Consider two
agents, whose information partitions are $\Part_1 = \{ \{ \omega, \omega
+ 1 \} \mid \omega \in X \} \cup \{ \{ \omega \} \mid \omega \notin X
\}$ and $\Part_2 = \{ \{ \omega \} \mid \omega < 1 \} \cup \{ [1,2]
\}$. Then $X \cup [1,2]$ is the block of the common-knowledge
partition that includes the block $[1,2]$ in $\Part_2$. From information
partitions composed of the simplest events---singletons, pairs, and a
closed interval---we have passed to a common-knowledge partition with
a block that, depending on what is $X$, might even be outside the
projective hierarchy. Knowledge of an event by individual agent is
likewise problematic. In the present example, the event that agent 1
knows $(0,1)$ is $(0,1) \setminus X$.

These measurability problems dictate that information partitions
should be represented as equivalence relations. If $E_1$ and $E_2$ are
$\BS^1_1$ (that is, analytic) equivalence relations, then the meet of
the partitions that they induce is induced by the transitive closure
of their union. This transitive closure is also a $\BS^1_1$
equivalence relation.\footnote{Composition is defined with a single
  existential quantifier, and thus takes a pair of $\BS^1_1$
  relations to a $\BS^1_1$ relation. The countable union of
  $\BS^1_1$ relations is $\BS^1_1$. Cf.~\citet[Theorem 2B.2,
    p.~54]{Moschovakis-2009}.}
Moreover, if an information partition is represented by a $\BS^1_1$
relation and an event is $\BP^1_1$ (that is, co-analytic), then
knowledge of the event according to the information partition is also
a $\BP^1_1$ event.\footnote{This is equivalent, by \eqn o, to the fact
  that the saturation of a $\BS^1_1$ set with respect to a $\BS^1_1$
  equivalence relation is $\BS^1_1$. This latter fact is true because
  the saturation is defined with a single extential
  quantifier.} 

This observation implies that knowledge of a Borel event is a
universally measurable event---surely a threshold condition for
incorporating the analysis of knowledge into Bayesian theory.
However, it is easy to envision a noncooperative Bayesian game
situation in which a player must condition on the event that some
other event is common knowledge.\footnote{In a game-theoretic analysis
  of the ``coordinated-attack problem'' (cf.~\citet{Rubinstein-1989}),
  there is an equilibrium in which two players will take complementary
  actions iff it is common knowledge that each has received a
  signal. The two signals are privately observed by the respective
  players and are highly, but imperfectly, correlated. Thus the
  players must exchange infinitely many messages to one another (that
  is, `I have received the signal', `I have received your confirmation
  of receipt of the signal', `I have received your confirmation of
  receipt of my confirmation of receipt of the signal',\dots) in order
  to attain common knowledge. If the communications channel closes at
  any stage, then common knowledge is not reached, and the action will
  not be taken. Closure of the channel is irreversible. The game can
  be elaborated by supposing that the stochastic distribution of
  payoffs from coordination is dependent on whether common knowledge
  has been achieved. (For example, the communication channel may have
  to be used, and therefore must not be broken, in order to accomplish
  the task. The probability that the channel being open at that time
  is positive if common knowledge has been achieved, but is zero
  otherwise.) Thus, in order for the game theorist to prove that
  coordination conditional on attaining common knowledge is an
  equilibrium, the expected payoff from coordination conditional on
  common knowledge being attained must be well defined.} Then the game
theorist needs to model the event that both players' receipt of their
respective signals is common knowledge as being a $\BD^1_1$ (that is,
Borel) event, not just a $\BP^1_1$ event. It might be thought possible to
avoid this difficulty by making tighter modeling assumptions regarding
both the agents' information structures and also the event to which
the common-knowledge operator is to be applied.

In particular, in most applications to Bayesian decision
theory and game theory, it is reasonable to specify each agent's
information as a $\BD^1_1$ equivalence relation, or even as a smooth
or closed Borel relation rather than as an arbitrary $\BS^1_1$
equivalence relation.\footnote{Smoothness (also called tameness) and
  closedness are co-extensive for equivalence relations on standard
  Borel spaces. Cf.~\citet[proof of Theorem 1.1,
    p.~920]{HarringtonKechrisLouveau-1990}. Standard Borel spaces are
  defined below, in \sec d.} Thus it may be asked: if
the graphs of $E_1$ and $E_2$ are in $\BD^1_1$ or in some smaller
class, then how is the graph of the transitive closure of $E_1 \cup
E_2$ restricted?

It will be shown here that no significant restriction of the
common-knowledge partition is implied by such restriction of agents'
information partitions. This finding is not surprising, since
restricting the complexity of individuals' equivalence relations does
not obviate the use of an existential quantifier to define the
transitive closure of a relation. Nevertheless, the syntactic form of
a specific description of a set does not determine the intrinsic
complexity of the set, so it needs to be shown that common-knowledge
equivalence relations derived from Borel equivalence relations are not
lower in the projective hierarchy, as a class, than their definition
would suggest. Moreover, \prp i will show that being the union
of finitely many (in fact, of fewer that $2^{\aleph_0}$) Borel
  equivalence relations---that is, representability in the form, of
  which the transitive closure is an equivalence
  relation specifying common-knowledge---is a stronger property than
  being an arbitrary Borel, reflexive, symmetric relation.

To define the transitive closure of $R \subseteq \Omega \times
\Omega$, let $R^{(1)} = R$ and $R^{(n+1)} = R \cmps R^{(n)}$ (that is,
the composition of relations $R$ and $R^{(n)}$). Denote the transitive
closure of $R$ by $\tc{R} = \bigcup_{n \in {\Nat_+}} R^{(n)}$.  It will be
proved here that, if $\Omega$ is a Polish space and $E_0 \subset
\Omega \times \Omega$ is a $\BS^1_1$ equivalence relation, then there
are smooth $\BD^1_1$ equivalence relations $E_1$ and $E_2$ and a
$\BD^1_1$ subset $Z$ of $\Omega$, such that $\tc{E_1 \cup E_2} \rst Z$
is Borel equivalent to $E_0$.\footnote{$R \rst Z = R \cap (Z \times
  Z)$. Let restriction take precedence over Boolean operations. For
  example, $X \cup R \rst Z \cap Y$ means $X \cup(R \rst Z) \cap X$.}
If $\Omega$ is the Baire space, then $E_1$ and $E_2$ can be taken to
be closed, $Z$ can be taken to be open, and the Borel equivalence can
be taken to be a homeomorphic equivalence.

\Section b{The case of the Baire space}

First take $\Omega$ to be the Baire space, $\N =
\Nat^\Nat$.\footnote{$\Nat = \{0,1,\ldots\}$. $\N$ is topologized as
  the product of discrete spaces.} Define subsets
$X$ and $Y$ of $\N$ by $X = \{ \alpha | \alpha_0>0 \}$ and $Y = \{
\alpha | \alpha_0 = 0 \}$.  $X$ and $Y$ are both homeomorphic to $\N$,
and homeomorphisms $f\from X \to Y$ and $g\from Y\times Y \times Y \to
Y$ are routine to construct.\footnote{Since $Y$ is homeomorphic with
  $\N$, $g$ can be constructed from the function described by
  \citet[p.~31]{Moschovakis-2009}.} Each of $X$ and $Y$ is both open
and closed in $\N$. It follows that, if $Z$ is either $X$ or $Y$, then
$A \subseteq Z$ is open (resp.~closed, Borel, $\BS^1_1$) as a subset
of $A$ iff it is open (resp.~closed, Borel, $\BS^1_1$) as a subset of
$Z$. This invariance to the ambient space extends to product
spaces. (For example a subset of $X \times Y$ is closed in $X \times
Y$ iff it is closed in $\N \times \N$.) In subsequent discussions,
subsets of these subspaces will be characterized (for example, as
being closed) without mentioning the subspace.

\theor g{If $E \subseteq X \times X$ is a $\BS^1_1$ equivalence
  relation, then there are equivalence relations $I$ and $J$ on $\N
  \times \N$, each of which has a closed graph, such that $E = \tc{I
    \cup J} \rst X$.}

Before proceeding to the proof of this theorem, note that $I \cup J$
is a closed, reflexive, symmetric relation. Thus, \thm g has the
following corollary.

\corol h{If $E \subseteq X \times X$ is a $\BS^1_1$ equivalence
  relation, then there is a closed, reflexive, symmetric relation $R$ on $\N
  \times \N$, such that $E = \tc{R} \rst X$.}

\uppercase\thm g would follow from \cor h if every closed,
reflexive, symmetric relation were the union of two closed equivalence
relations, but that is not the case.

\propo i{Let $\alpha \in \N$. Define $R = D \cup
  (\{ \alpha \} \times \N) \cup (\N \times \{ \alpha \})$, and define 
\[ \mathcal{E} = \bigcup_{} \left\{ (D \cup 
\{ (\alpha, \beta), (\beta, \alpha) \}) | \beta \in \N \setminus
\{ \alpha \} \right\}. \]
$R = \bigcup \mathcal{E}$; every $E \in \mathcal{E}$ is an
equivalence relation; $R$ is closed, reflexive, and symmetric; and
$2^{\aleph_0}$ is the cardinality of $\mathcal{E}$. There is no other
set $\mathcal{F}$ of equivalence relations such that $R =
\bigcup\mathcal{F}$. Thus, $R$ is not a union of fewer that
$2^{\aleph_0}$ equivalence relations.}

\begin{proof} 
The assertions regarding $\mathcal{E}$ are obvious from
its construction. To obtain a contradiction from supposing that
$\mathcal{E}$ were not unique, suppose that $R$ were also the union of
a set $\mathcal{F} \neq \mathcal{E}$ of equivalence relations. Not
$\mathcal{F} \subsetneq \mathcal{E}$. So, there must be some $E \in
\mathcal{F} \setminus \mathcal{E}$. By symmetry, there must be three
distinct points, $\alpha, \beta, \gamma$ such that $\{ (\beta,
\alpha), (\alpha, \gamma) \} \subseteq E$. Since $E$ is transitive,
$(\beta, \gamma) \in E \setminus R$, contrary to $R = \bigcup
\mathcal{F}$.
\end{proof}

\Section c{Proof of the theorem}

Denote the diagonal (that is, identity) relation in $\N \times \N$ by $D =
\{ (\alpha, \alpha) | \alpha \in \N \}$. $D$ is closed.

If $1 \le i < j \le k$ and $\vec \alpha = (\alpha_1,
\alpha_2,\ldots,\alpha_k) \in \N^k$, then a transposition mapping is
defined by $\tpos[ij]{\vec \alpha} = (\alpha_1, \ldots,\alpha_{i-1},
\alpha_j, \alpha_{i+1},\ldots, \alpha_{j-1}, \alpha_i,
\alpha_{j+1},\ldots, \alpha_k)$.\footnote{A sub-sequence of
  subscripted alphas distinct from $\alpha_i$ and $\alpha_j$ having
  subscripts that are not increasing, which occurs if $i=1$ or $j =
  i+1$ or $j=k$, denotes the empty sequence.} The abbreviation $\Tpos
A = \tpos{A} = \{ \tpos{\alpha} | \alpha \in A \}$ will sometimes be
used. Each $\transpose[ij]$ is a homeomorphism of $\N^k$ with
itself. Note that $\transpose[ij] \rst X$ and $\transpose[ij]
\rst Y$ map $X^k$ and $Y^k$ homeomorphically onto themselves.

Recall that a relation $E \subseteq X \times X$ is $\mathbf{\Sigma}^1_1$
iff there is a set $F$ such that
\display j{F \subseteq X \times X \times \N \text{ is closed, and }
(\alpha, \beta) \in E \iff \exists \gamma \; (\alpha,
  \beta,\gamma) \in F.}
\lema h{If $E \subseteq X \times X$ is symmetric, then $E$ is
$\mathbf{\Sigma}^1_1$ iff there is a closed, $\transpose$-invariant
set $F \subset X \times X \times X$ that satisfies \eqn j.}

\begin{proof}
Let $F_0$ satisfy \eqn j. Let $h$ be a homeomorphism from $\N$ to $X$,
and  define $F_1 \subseteq X \times X \times X$ by $(\alpha, \beta,
\gamma) \in F_0 \iff (\alpha, \beta, h(\gamma)) \in F_1$. $F_1$ also
satisfies \eqn j, then, and it is closed. By symmetry of $E$, 
$\Tpos{F_1}$ is another closed set that satisfies \eqn
j. Consequently, $F = F_1 \cup \Tpos{F_1}$ is a $\transpose$-invariant
closed set that satisfies \eqn j.
\end{proof}

Let $\id$ denote the identity function on $\N$. If $K,L,M,N$ are any
sets, and $p\from  K \to L$ and $q\from  M \to N$, then denote the
product mapping by $p \times q \from  K \times M \to L \times N$.

The two closed equivalence relations that \thm g asserts to exist are
defined from the homeomorphisms $f$ and $g$ introduced in \sec b, and
the closed, $\transpose$-invariant set $F$ guaranteed to exist by \lem
h, as follows.
\display l{\begin{split}
j(\alpha, \beta, \gamma) &= g(f(\alpha), f(\beta), f(\gamma)).\\
G &= \{ (\alpha, j(\alpha, \beta, \gamma)) | (\alpha,
  \beta, \gamma) \in F \} \subseteq X \times Y;\\
H &= \{ (j(\alpha, \beta, \gamma), j(\beta, \alpha, \gamma)) |
  (\alpha, \beta, \gamma) \in X \times X \times X \};\\ 
I &= D \cup G \cup \Tpos{G} \cup \Tpos{G} \cmps G;\\
J &= D \cup H.
\end{split}}

\lema j{$D$, $G$, $\Tpos{G}$, $H$, and $J$ are closed.}

\begin{proof}
$D$ is closed because $\N$ is a metric space.

The function $i \times j$ is a homeomorphism from $X \times X \times X
\times X$ to $X \times Y$. Being a homeomorphism, it is an open
mapping (which takes closed sets to closed sets). $G = [i \times
  j](((D \rst X \times X \times X)) \cap (X \times F))$.  $D \rst
X \times X \times X$ and $X \times F$ are both closed subsets of $X
\times X \times X \times X$, so $G$ is closed. $\Tpos{G}$ is closed,
as the image of $G$ under $\transpose$, a self-homeomorphism of $\N
\times \N$.

$j \times j$ is a homeomorphism from $(X \times X \times X) \times (X
\times X \times X)$ to $Y \times Y$. The image under $j \times j$ of a
closed subset of its domain is therefore closed in its
range. $\{ ((\alpha, \beta, \gamma), (\beta, \alpha, \gamma)) |
  (\alpha, \beta, \gamma) \in X \times X \times X \}$ is $\transpose[23]
\circ \transpose[25] (D \rst X
\times D \rst X \times D \rst X)$, which is
closed. $H$, the image of this set under $j \times j$, is therefore
closed.

$J$, the union of two closed sets, is closed.
\end{proof}

\lema k{$G \cmps \Tpos{G} = D \rst X$.  $\Tpos{G} \cmps G = \{
  (j(\alpha, \beta, \gamma), j(\alpha, \delta, \epsilon)) | (\alpha,
  \beta, \gamma) \in F \text{ and } (\alpha, \delta, \epsilon) \in F
  \}$. $H = \Tpos{H}$. $H^{(2)} = D \rst Y$.  $G \cmps H = \{ (\alpha,
  j(\beta, \alpha, \gamma) | (\alpha, \beta, \gamma) \in F \}$. $G
  \cmps H \cmps \Tpos{G} = E$.  $\Tpos{G} \cmps G$ and $I$ are
  closed.}

\begin{proof}
All assertions except the one regarding closedness of $\Tpos{G}
\cmps G$ and $I$ are verified by straightforward calculations. That $F$ is
invariant under $\transpose$ is used to show that $H = \Tpos{H}$ and
that $E \subseteq G \cmps H \cmps \Tpos{G}$.

The proof that $\Tpos{G} \cmps G$ is closed is parallel to the proof that
$H$ is closed. According to the first part of this lemma, 
$\Tpos{G} \cmps G = [j \times j](\transpose[24](D \rst X
\times X \times X \times X \times X) \cap (F \times F))$.

$I$, the union of four closed sets, is closed.
\end{proof}

\lema l{$I$ and $J$ are equivalence relations.}

\begin{proof}
These relations are reflexive and symmetric, so their transitive
closures are equivalence relations. Thus, the lemma is equivalent to
the assertion that $I = \tc{I}$ and $J = \tc{J}$. For any relation
$K$, $K^{(2)} = K$ is sufficient for $K = \tc{K}$. In the following
calculations of $I^{(2)}$ and $J^{(2)}$, composition of relations is
distributed over unions. Terms that evaluate by
identities that were calculated in \lem k to a previous term or its
sub-relation, are omitted from the expansion by terms in the
pentultimate step of each calculation.
\display m{\begin{split}
I^{(2)} &= (D \cup G \cup \Tpos{G} \cup \Tpos{G} \cmps G) \cmps
(D \cup G \cup \Tpos{G} \cup \Tpos{G} \cmps G)\\
&= (D \cup G \cup \Tpos{G} \cup \Tpos{G} \cmps G) \cup (G \cup G \cmps
\Tpos{G} \cup G \cmps \Tpos{G} \cmps G) \cup (\Tpos{G} \cup \Tpos{G}
\cmps G \cup \Tpos{G} \cmps \Tpos{G} \cup \cup \Tpos{G} \cmps \Tpos{G}
\cmps G)\\
&\qquad \cup (\Tpos{G} \cmps G \cup \Tpos{G} \cmps G \cmps G \cup
\Tpos{G} \cmps G \cmps \Tpos{G} \cup \Tpos{G} \cmps G \cmps \Tpos{G}
\cmps G)\\
&= D \cup G \cup \Tpos{G} \cup \Tpos{G} \cmps G\\
&= I.\\
&\strut\\
J^{(2)} &= (D \cup H)(D \cup H)\\
&= (D \cup H) \cup (H \cup H^{(2)})\\
&= D \cup H\\
&= J.
\end{split}}
\end{proof}

\begin{proof}[Proof of \thm g]

Lemmas \ref{lem:j}--\ref{lem:l} show that the each of the relations
$I$ and $J$ on $\N \times \N$, is an equivalence relation that has a
closed graph. It remains to be shown that that $E = \tc{I \cup J} \cap
(X \times X)$. Note that, since $D \subseteq I \cup J$, $I \cup J
\subseteq (I \cup J)^{(2)} \subseteq (I \cup J)^{(3)} \subseteq \ldots\;$
Hence, if $(I \cup J)^{(n)} = (I \cup J)^{(n+1)}$, then $(I
\cup J)^{(n)} = \tc{I \cup J}$.

The following calculation shows that $(I \cup J)^{(5)} = (I \cup
J)^{(6)}$. The calculation is done recursively, according to the
following recipe at each stage $n>1$:
\begin{enumerate}
\item
Begin with the equation $(I \cup J)^{(n+1)} = (I \cup J)(I \cup J)^{(n)}$.
\item
Rewrite $(I \cup J)$ as $D \cup G \cup \Tpos{G} \cup \Tpos{G}G \cup H$
according to \eqn l, rewrite $(I \cup J)^{(n)}$ according to the
result of the previous step, and then distribute composition of
relations over union in the resulting equation.
\item
For each identity stated in \lem k, and for each identity that, for
some $K \in \{G, \Tpos{G}, H \}$, equates a composition $K \cmps D$ or
$D \cmps K$ of $K$ and $D$ (or a restriction of $D$ to a product set
of which $K$ is a subset) to $K$, do as follows: Going from left to
right, apply the identity wherever possible.\footnote{Let $P=D \rst X$
  and $Q=D \rst Y$. Identities are applied in the following order at
  each stage of the recursion: $D \cmps D = D$, $D \cmps E = E$, $D
  \cmps G = G$, $D \cmps \Tpos{G} = \Tpos{G}$, $D \cmps H = H$, $D
  \cmps P = P$, $D \cmps Q = Q$, $E \cmps D = E$, $E \cmps E = E$, $E
  \cmps P = E$, $G \cmps D = G$, $G \cmps \Tpos{G} = P$, $G \cmps H
  \cmps \Tpos{G} = E$, $G \cmps Q = G$, $\Tpos{G} \cmps D = \Tpos{G}$,
  $\Tpos{G} \cmps P = \Tpos{G}$, $H \cmps D = H$, $H \cmps H = Q$, $H
  \cmps Q = H$, $P \cmps D = P$, $P \cmps E = E$, $P \cmps G = G$, $P
  \cmps P = P$, $Q \cmps D = Q$, $Q \cmps \Tpos{G} = \Tpos{G}$, $Q
  \cmps H = H$, $Q \cmps Q = Q$.}  Repeat this entire step (consisting
of one pass per identity) until no further simplifications are
possible.
\item
Delete compositions of relations that include terms $K \cmps L$ such
that the range of $K$ and the domain of $L$ (viewed as
correspondences) are disjoint, in which case the term denotes the empty
relation. Delete $D \rst X$ (occurring as a term by
itself), of which $D$ is a superset. 
\item
Delete each term of form $[K] \cmps \Tpos{G} \cmps [L]$
(resp.~$[K] \cmps G \cmps [L]$) from a union in which the
corresponding term for its superset, $[K] \cmps \Tpos{G} \cmps E \cmps
[L]$ (resp.~$[K] \cmps E \cmps G \cmps [L]$) also appears.
(One or
both of the bracketed sub-terms may be absent from both terms in the
pair.) Delete $D$ (occurring as a term by itself) from every union
that contains both $D \rst Y$ and $E$, since $D \subseteq D
\rst Y \cup E$.
\item
Reorder terms lexicographically, in the order $D< D \rst Y
<E<G<\Tpos{G}<H$. Delete repeated terms.
\end{enumerate}

\display n{\begin{split} (I \cup J)^{\hphantom{(2)}} &= D \cup G \cup
    \Tpos{G} \cup \Tpos{G}G \cup H\\ 
&\strut\\
(I \cup J)^{(2)} &= D \cup D \rst Y \cup G \cup GH \cup
    \Tpos{G} \cup \Tpos{G}G \cup \Tpos{G}GH\\
&\cup H \cup H\Tpos{G} \cup H\Tpos{G}G\\  
&\strut\\ 
(I \cup J)^{(3)} &= D \rst Y \cup E \cup EG \cup GH \cup
    \Tpos{G}E \cup \Tpos{G}EG \cup \Tpos{G}GH\\ 
&\cup H \cup H\Tpos{G} \cup H\Tpos{G}G \cup H\Tpos{G}GH\\
&\strut\\
(I \cup J)^{(4)} &= D \rst Y \cup E \cup EG \cup EGH \cup
    \Tpos{G}E \cup \Tpos{G}EG \cup \Tpos{G}EGH\\ 
&\cup H \cup H\Tpos{G}E \cup H\Tpos{G}EG \cup H\Tpos{G}GH\\
&\strut\\
(I \cup J)^{(5)} &= D \rst Y \cup E \cup EG \cup EGH \cup
    \Tpos{G}E \cup \Tpos{G}EG \cup \Tpos{G}EGH\\
&\cup H \cup H\Tpos{G}E \cup H\Tpos{G}EG \cup H\Tpos{G}EGH\\
&= (I \cup J)^{(6)}
\end{split}}

Thus $\tc{I \cup J} = (I \cup J)^{(5)}$. Note that $D \rst Y$, $G$,
$\Tpos{G}$, $H$ and all relations of form  or 
$\Tpos{G} \cmps Q$ or $H \cmps Q$ or $Q \cmps G$ or $Q \cmps H$ (where
variable $Q$ ranges over compositions of $G$, $\Tpos{G}$, $H$, and
$E$), are disjoint from $X \times X$. Therefore, from the calculation
in \eqn n of $(I \cup J)^{(5)}$ as a union of $E$ with such relations,
it follows that $\tc{I \cup J} \cap (X \times X) = E$.
\end{proof}

\Section d{The general case of a standard Borel space}

In this concluding section, \thm g is generalized in two ways to an
arbitrary standard Borel space. A \emph{standard Borel space} is a
pair $\BO_0 = (\Omega_0, \mathcal{B}_0)$ such that, for some pair $\BO
= (\Omega, \mathcal{B})$, $\mathcal{B}$ is the $\sigma\/$-algebra of
Borel subsets of the set $\Omega$ under some Polish topology,
$\Omega_0 \in \mathcal{B}$, and $\mathcal{B}_0 = \{ B_0 | \, \exists B \,
[B \in \mathcal{B} \text{ and } B_0 = B \cap \Omega_0] \}$. A
\emph{Borel isomorphism} of standard 
Borel spaces $\BO_0$ and $\BO$ is a $\BD^1_1$ function $k \from
\Omega_0 \to \Omega$ such that $k^{-1} \from \Omega \to \Omega_0$
exists and is also $\BD^1_1$. A $\BD^1_1$ subset of a standard Borel
space is also a standard Borel space, and every two uncountable
standard Borel spaces are
isomorphic.\footnote{\citet[pp.~338--9]{Mackey-1957}. Henceforth,
  $\mathcal{B}$ will be implicit and the structure $\BO$ will be
  identified with the set $\Omega$ on which it is defined.}

In both generalizations, the concept of smoothness of a Borel
equivalence relation substitutes for the concept of closedness that
appears in \thm g. If $E \subseteq \Omega
\times \Omega$ is a $\BD^1_1$ equivalence relation, and if there is a
set $\{ Y_n \}_{n \in \Nat}$ of $\BD^1_1$ sets such that $(\omega,
\omega') \in E \iff \forall n \; [\omega \in Y_n \iff \omega' \in
  Y_n]$, then $E$ is a \emph{smooth} equivalence relation. By
\citet[proof of Theorem 1.1, p.~920]{HarringtonKechrisLouveau-1990},
every equivalence relation with closed graph is smooth. If $k \from
\Omega_0 \to \Omega$ is $\BD^1_1$ and $E \subseteq \Omega \times
\Omega$ is a smooth $\BD^1_1$ equivalence relation, then $E_0 \subseteq
\Omega_0 \times \Omega_0$ defined by $(\psi, \omega) \in E_0 \iff (k(\psi),
k(\omega)) \in E$ is also smooth, with $E_0\/$-equivalence determined by
$\{k^{-1}(Y_n) \}_{n \in \Nat}$.

The first generalization of \thm g asserts Borel embeddability of
an arbitrary $\BS^1_1$ equivalence relation. If $\Omega_0$ and $\Omega$
are standard Borel spaces, and $E_0 \subseteq \Omega_0 \times \Omega_0$ and
$E \subseteq \Omega \times \Omega$ are $\BS^1_1$ equivalence
relations, then a \emph{Borel embedding} of $E_0$ into $E$ is a Borel
isomorphism $e \from \Omega_0 \to Z \subseteq \Omega$ that extends
naturally to a Borel isomorphism from $E_0$ to $E \rst Z$. That is,
$(\psi, \omega) \in E_0 \iff (e(\psi), e(\omega)) \in E$.

\corol m{Let $\Omega_0$ and $\Omega$ be  standard Borel spaces, and let $E_0
  \subseteq \Omega_0 \times \Omega_0$ be a $\BS^1_1$ equivalence
  relation. There are smooth $\BD^1_1$ equivalence relations $E_1 \subseteq
  \Omega \times \Omega$ and $E_2 \subseteq \Omega \times \Omega$ such
  that $E_0$ is Borel embeddable in $\tc{E_1 \cup E_2}$.}

\begin{proof}
If $\Omega_0$ is countable, then $E_1$ and $E_2$ can both be taken to
be the image of $E_0$ under an arbitrary injection of $\Omega_0$ into
$\Omega$. Otherwise, there is a Borel isomorphism $k_0 \from \Omega_0
\to X$ (where $X$ is as in \thm g), and there is a Borel isomorphism
$k \from \Omega \to X$. Define $e = k^{-1} \circ k_0$ and define $Z
\subseteq \Omega$ by $Z = e(\Omega_0)$. If $E \subset X \times X$ is
defined by $(\alpha, \beta) \in E \iff (k^{-1}_0(\alpha),
k^{-1}_0(\beta)) \in E_0$, then $E$ is a $\BS^1_1$ equivalence
relation.\footnote{\citet[Theorem 2B.2, p.~54]{Moschovakis-2009}.} Let
$I$ and $J$ be the closed equivalence relations defined in \eqn l, and
define $(\psi, \omega) \in E_1 \iff (k(\psi), k(\omega)) \in I$
and $(\psi, \omega) \in E_2 \iff (k(\psi), k(\omega)) \in
J$. $E_1$ and $E_2$ are smooth. Now the corollary follows immediately
from \thm g.
\end{proof}

The second generalization of \thm g applies to a $\BS^1_1$ equivalence
relation that, in a sense, does not occupy the entire product space
$\Omega \times \Omega$. Specifically, the set of points, the
singletons of which are blocks of the partition induced by the
relation, must have an uncountable $\BD^1_1$ subset.\footnote{If $E$
  is a $\BS^1_1$ equivalence relation, then the set of all such points
  is a $\BP^1_1$ subset of $\Omega$. One sufficient condition for an
  uncountable $\BP^1_1$ set, $W$, to have an uncountable $\BD^1_1$
  subset is that there should be a nonatomic measure, $\mu$, on
  $\Omega$ such that $\mu^*(\Omega \setminus W) < \mu(\Omega)$ (where
  $\mu^*$ is outer measure). Another sufficient condition is that $W$
  should have a perfect (hence both uncountable and $\BD^1_1$) subset.
  Two sufficient conditions for every uncountable $\BP^1_1$ set to
  have a non-empty perfect subset---albeit conditions that are
  independent of ZFC set theory (if ZFC is consistent)---are provided
  by \citet[Exercise 6G.10, p.~288; and Corollary 8G.4,
    p.~419]{Moschovakis-2009}. (\citet[Theorem
    25.38. p.~499]{Jech-2002} states the ``boldface'' implication of
  Corollary 8G.4.) It is provable in ZFL that there is an uncountable
  $\BP^1_1$ set (in fact, a $\Pi^1_1$ set) without a non-empty
  perfect subset. \citet[Exercise 5A.8, p.~212]{Moschovakis-2009}.}

\corol n{Suppose $\Omega$ is a standard Borel space and that $E
  \subseteq \Omega \times \Omega$ is a $\BS^1_1$ equivalence relation
  such that, for some uncountable $\BD^1_1$ set $B \subseteq \Omega$,
  $E \rst B = D \rst B$. Define $\Omega_0 = \Omega \setminus B$.  Then
  there are smooth $\BD^1_1$ relations $E_1$ and $E_2$, such that
  $E \setminus D \subseteq E \rst \Omega_0 \cup D \rst B$}

Finally, \cor n provides a negative answer to the question, raised in
the introduction, of whether the saturations of Borel sets (or even of
singletons) with respect to the transitive closures of unions of
smooth Borel equivalence relations lie within any significantly
restricted sub-class of $\BS^1_1$.

\corol o{Suppose $\Omega$ is a standard Borel space and that $S
  \subseteq \Omega$ is a $\BS^1_1$ set such that, for some $\BD^1_1$
  set $\Omega_0$, $S \subseteq \Omega_0$ and $\Omega \setminus
  \Omega_0$ is uncountable. Then there are smooth $\BD^1_1$ relations
  $E_1$ and $E_2$, such that for every non-empty $A \subseteq S$,
  $[A]_{\tc{E_1 \cup E_2}} \cap \Omega_0 = S$.}

\begin{proof}
Define $(\psi, \omega) \in E \iff [ \{ \psi, \omega \} \subseteq S
  \text{ or } \psi = \omega]$, specify $B = \Omega \setminus
\Omega_0$, and apply \cor n.  \relax For some block, $\pi$, of the
partition induced by $\tc{E_1 \cup E_2}$, $\pi \cap \Omega_0 = S$. Therefore,
if $\emptyset \neq A \subseteq S$, then $[A]_{\tc{E_1 \cup E_2}} \cap \Omega_0 = S$.
\end{proof}

\bibliographystyle{plainnat}

\begin{thebibliography}{8}
\providecommand{\natexlab}[1]{#1}
\providecommand{\url}[1]{\texttt{#1}}
\expandafter\ifx\csname urlstyle\endcsname\relax
  \providecommand{\doi}[1]{doi: #1}\else
  \providecommand{\doi}{doi: \begingroup \urlstyle{rm}\Url}\fi

\bibitem[Aumann(1976)]{Aumann-1976}
Robert~J. Aumann.
\newblock Agreeing to disagree.
\newblock \emph{Ann.~Statist.}, 4\penalty0 (6):\penalty0 1236--1239, 1976.
\newblock URL \url{http://projecteuclid.org/euclid.aos/1176343654}.

\bibitem[Harrington et~al.(1990)Harrington, Kechris, and
  Louveau]{HarringtonKechrisLouveau-1990}
L.~A. Harrington, A.~S. Kechris, and A.~Louveau.
\newblock A \uppercase{G}limm-\uppercase{E}ffros dichotomy for
  \uppercase{B}orel equivalence relations.
\newblock \emph{J.~Amer.~Math.~Soc.}, 3\penalty0 (4):\penalty0 903--928, 1990.

\bibitem[Jech(2002)]{Jech-2002}
Thomas Jech.
\newblock \emph{Set Theory}.
\newblock Springer-Verlag, third edition, 2002.

\bibitem[Kripke(1959)]{Kripke-1959}
Saul~A. Kripke.
\newblock A completeness theorem in modal logic.
\newblock \emph{J.~Symbolic Logic}, 24\penalty0 (1):\penalty0 1--14, 1959.

\bibitem[Mackey(1957)]{Mackey-1957}
G.~W. Mackey.
\newblock Borel structure in groups and their duals.
\newblock \emph{Trans. Amer. Math. Soc.}, 85\penalty0 (1):\penalty0 134--165,
  1957. \hfil
\newblock URL \texttt{http://www.ams.org/journals/tran/1957-085-01/\hfil\break S0002-9947-1957-0089999-2/home.html}.

\bibitem[Moschovakis(2009)]{Moschovakis-2009}
Yannis~N. Moschovakis.
\newblock \emph{Descriptive Set Theory}.
\newblock American Mathematical Association, second edition, 2009.
\newblock URL
  \url{http://www.math.ucla.edu/\textasciitilde ynm/books.htm}.

\bibitem[Pohlers(2009)]{Pohlers-2009}
Wolfram Pohlers.
\newblock \emph{Proof Theory: The First Step into Impredicativity}.
\newblock Springer-Verlag, 2009.

\bibitem[Rubinstein(1989)]{Rubinstein-1989}
Ariel Rubinstein.
\newblock The electronic mail game: strategic behavior under ``almost common
  knowledge''.
\newblock \emph{Amer. Econ. Rev.}, 79\penalty0 (3):\penalty0 385--391, 1989.

\end{thebibliography}

\end{document}